\documentclass[a4paper,10pt]{article}

\usepackage{amsmath,amssymb} 
\usepackage{amsthm}        

\newcommand{\N}{\mathbb{N}}

\newcommand{\F}{\mathbb{F}}
\newcommand{\C}{\mathbb{C}}

\newcommand{\spec}{\mathop{\rm Spec}\nolimits}
\newcommand{\ind}[2]{#1_{\scriptscriptstyle{#2}}}
\newcommand{\parenth}[1]{\left(#1\right)}
\newcommand{\croch}[1]{\left[#1\right]}
\newcommand{\card}{\mathop{\rm Card}\nolimits}
\newcommand{\cjc}{\mathop{\rm CJC}\nolimits}
\newcommand{\njc}{\mathop{\rm NJC}\nolimits}

\title{On isomorphisms of factorial domains and the Jacobian conjecture in any characteristic} 
\author{
Kossivi \textsc{Adjamagbo}\\
\footnotesize{Universit\'e de Paris VI, UFR 929}\\
\footnotesize{Institut de Math\'ematiques de Jussieu, UMR 7586 CNRS}\\
\footnotesize{4, place Jussieu, 75252 PARIS Cedex 05}\\
\footnotesize{e-mail : adja@math.jussieu.fr}
}
\date{} 

\begin{document}

\maketitle

\begin{abstract}
The main theorem (2.2) consists in two characterizations of isomorphisms of factorial domains in terms of prime or primary rings elements, and unramified, flat or weakly injective affine schemes morphisms. In order to apply this theorem to the famous Jacobian Conjecture, we first introduce its different versions in any characteristic (3.1), and give two reformulations of some these versions in terms of domains of positive characteristic (3.8) and finite prime fields (3.9). Finally, we deduce from the main theorem an original reformulation of the any characteristic version of the Jacobian Conjecture in terms of prime or primary rings elements (3.11).
\end{abstract}

\section*{Introduction}

The main aim of the present work is to  establih a ``U.F.D. isomorphism theorem" (see 2.2 below), which gives two characterizations of factorial domains isomorphisms in terms of prime or primary elements an unramified, flat or ``weakly injective" schemes morphism. More precisely, we prove that a morphism $\phi$ between U.F.D.'s of characteristic $p$ $A$ and $B$ is an isomorphism if and only if $\phi$ is a separable morphism of finite type such that its co-morphism $\phi^{*}$ from $\spec(B)$ to $\spec(A)$ is ``weakly injective" (\emph{i.e.} injective on the inverse image by $\phi^{*}$ of the set of invertible elements of the ring $R$) and $[B_*^{-1}B : \phi(A_*)^{-1}\phi(A)]\not\in p\N$ (where we write $R_*$ for $R-\{0\}$ for any ring $R$ or $R=\N$), or if and only if $\phi$ is flat morphism (hence a monomorphism) transforming prime elements of $A$ into prime elements of $B$ such that $B^*=\phi(A^*)$ and $B_*^{-1}B$ is a finite and separable field of $\phi(A_*)^{-1}\phi(A)$ with a degree not divisible by $p$.

For clarity and to keep the exposition of the main theorem self-contained, in part 1 of this paper, we recall basic definitions and results about separable algebras we need, in particular the ``first isomorphism theorem for algebras over a U.F.D." of \cite{AD1} (see 1.6 below), which is a jacobian and monogenic characterization of isomorphisc extensions of an U.F.D. and the basic tool of the proof.

In order to apply fully this main theorem to the Jacobian Conjecture in any characteristic, in the last part of this paper, we begin to precise various versions of the ``classical Jacobian Conjecture in any characteristic" introduced in \cite{AD1},3.1 (see 3.1 below) and to establish two reformulations of the classical Jacobian Conjecture in characteristic zero in terms of positive characteristic or finite fields (see 3.8 and 3.9 below), thanks to the fundamental theorem of model theory proved by L\"os in 1955 (see the proof of 3.8) and an appropriate version of Bezout theorem published by S. Abhyankar in 1966, but ignored by specialists of intersection theory and ``rediscovered" and generalized by W. Fulton in 1980 (see the proof of 3.4 below). After this clarification, we exploit the factoriality of polynomial algebras over a U.F.D. to deduce from the main theorem a ``prime refomulation theorem" for Classical Jacobian Conjecture in any characteristic, in terms of prime or primary elements (see 3.11 below).

Although the contents of the present paper were ready a long time ago, at the period of preparation of \cite{AD1}, it has been written during the visit of the author in June 1995 at the University of Coimbra, in Portugal and was available as preprint of the Institut de Math\'ematiques de Jussieu of University Paris VI since October 1996 \cite{AD2}. So before entering the subject, I would like to express my gratitude to professors J. Vaillant and J. Carvalho e Silva for this pleasant and useful opportunity.

\section{Recall on separable algebras} 

\subsection{Generalities on separable algebras} (see for instance \cite{AD1}, 1.1)

For any ring $A$, let us recall that we write $A^*$ (resp. $A_*$) for the set of invertible (resp. non zero) elements of $A$, and similarly $\N_*$ for $\N-\{0\}$.

Let $\phi$ be a morphism from a commutative ring $A$ to another one $B$, $\pi$ the canonical map from $B\ind{\otimes}{A}B$ to $B$ induced by $\phi$ (\emph{i.e.} $\pi (b\otimes b)=bb'$ for all $(b,b')\in B^2$), and $\ind{\Omega}{B/A}$ the $B$-module of $A$-differentials (which is by definition $\ker\pi/(\ker\pi)^2$)

1) One says that $B$ is separable $A$-algebra, or that $B$ is separable over $A$ or that $\phi$ is separable, if $B\ind{\otimes}{A}B$-module $B$ induced by $\pi$ is projective.
	
2) If $\ker\pi$ is a finitely generated ideal of $B\ind{\otimes}{A}B$, the separability of $B$ over $A$ is equivalent to one of the following conditions :
	
\begin{enumerate}
	\item $\ind{\Omega}{B/A}=0$
	\item The ideal $\ker\pi$ of $B\ind{\otimes}{A}B$ is generated by an idempotent element
	\item For each prime ideal $p$ of $B$, $\ind{B}{p}$ is separable over $\ind{A}{\phi^{-1}(p)}$
	\item For each maximal ideal $m$ of $B$, $\ind{B}{m}$ is separable over $\ind{A}{\phi^{-1}(m)}$
	\item For each prime ideal $q$ in $A$, $\ind{B}{q}$ is separable over $\ind{A}{q}$
	\item For each maximal ideal $n$ in $A$ , $\ind{B}{n}$ is separable over $\ind{A}{n}$
	\item For each prime $q$ in $A$ with residue field $K(q)$, $B\ind{\otimes}{A}K(q)$ is separable over $K(q)$
\end{enumerate}
 
3) One says that $B$ is an unramified $A$-algebra, or that $B$ is unramified over $A$, or that $\phi$ is unramified, its co-morphism $\phi^{-1}$ from $\spec(B)$ to $\spec(A)$ is unramified, if $B$ is separable over $A$ and if $\ker\pi$ is a finitely generated ideal of $B\ind{\otimes}{A}B$
 
4) This condition on $\ker\pi$ is satisfied if $B$ is an $A$-algebra of essentially finite type, \emph{i.e.} a localization of an $A$-algebra of finite type.
 
5) If $S\subset A$ and $T\subset B$ are multiplicative sets such that $\phi(S)\subset T$, $\ind{\phi}{S}$ the map from $S^{-1}A$ to $T^{-1}B$ induced by $\phi$, $\ind{\pi}{S}$ the canonical map from $T^{-1}B\ind{\otimes}{S^{-1}A}T^{-1}B$ to $T^{-1}B$ induced by $\ind{\phi}{S}$, and if $\phi$ is separable (resp. $\ker\pi$ is a finitely generated ideal of $B\ind{\otimes}{A}B$), then $\ind{\phi}{S}$ is separable (resp. $\ker\ind{\pi}{S}$ is finitely generated ideal of $T^{-1}B\ind{\otimes}{S^{-1}A}T^{-1}B$)
 
6) If $A$ is a field, then the separability of $B$ over means tha $\ind{\dim}{A}B$ is finite and the ring $B\ind{\otimes}{A}L$ is reduced for any field extension $L$ of $A$, which also means that the $A$-algebra $B$ is isomorphic the product of a finite number of $A$-algebra of the form $A[T]/PA[T]$, where $P$ is a prime element, with a non-zero discriminant, of the $A$-algebra of polynomials $A[T]$ in one indeterminate $T$ over $A$.
 
7) If $A$ and $B$ are local rings, with respective maximal ideals $m(A)$ and $m(B)$, such that $\phi$ is local (\emph{i.e.} $\phi(m(A))\subset \phi(m(B))$) and $\ker\pi$ is a finitely generated ideal of $B\ind{\otimes}{A}B$, then $B$  is separable over $A$ if and if only $\phi(m(A))= \phi(m(B))$ and the field $B/m(B)$ is finite and separable over the field $A/m(A)$
  
8) Let $A[T]$ be the $A$-algebra of polynomials in one indeterminate $T$ over $A$ and $P$ an element of $A[T]\,\backslash\, A$ such that $B=A[T]/PA[T]$. $B$ is separable over $A$ if the discriminant  of $P$ is inversible in $A$. the converse is true if in addition $P$ in monic.
  
9) If $C$ is an commutative $A$-algebra and $B$ separable over $A$, then $B\ind{\otimes}{A}C$ is separable over $C$.

\subsection{Theorem (the first Jacobian Criterion of separability for a finitely presented algebra, see for instance \cite{WA}, th.5 or \cite{WR}, Prop.1.9)}

Let $A[X]$ be the algebra of polynomials generated by a system $X=\parenth{X_1,\ldots,X_n}$ of indeterminates over a commutative ring $A$, $P=\parenth{P_1,\ldots,P_m}\in A[X]^m$, $B$ the $A$-algebra $A[X]/\sum_{1\le i\le m}P_i A[X]$, and $J_X(P)=\parenth{\partial P_i/\partial X_j}\in M_{m,n}\parenth{A[X]}$ the jacobian matrix of $P$ with respect to $X$.

$B$ is separable over $A$ if and only if $m\ge n$ and $A[X]$ is generated as ideal by elements of $P$ and the maximal minors of $J_X(P)$.

\subsection{Theorem (the second Jacobian Criterion of separability for two finitely presented algebras, see \cite{AD1}, th.1.2).}

Let $K$ be a commutative ring, $K[Y]$ (resp.$K[X]$) the $K$-algebra of polynomials generated by a system $Y=\parenth{Y_1,\ldots,Y_n}$ (resp. $X=\parenth{X_1,\ldots,X_m}$) of indeterminates over $K$, $Q_1,\ldots,Q_s$ (resp. $P_1,\ldots,P_r$) elements of $K[Y]$ ($K[X]$), $\phi$ a $K$-algebra homomorphism from $A=K[Y]/\sum_{1\le i\le s}Q_iK[Y]$ to $B=K[X]/\sum_{1\le i\le r}P_iK[X]$ such that $\phi\parenth{Y_j+\sum_{1\le i\le s}Q_iK[Y]}=F_j+\sum_{1\le i\le r}P_iK[X]$, with $F_j\in K[X]$ for $1\le j\le n$, and $J_X(F,P)\in M_{n+r,m}(K[X])$ the jacobian matrix of $(F_1,\ldots,F_n,P_1,\ldots,P_r)$ with respect to $X$. 

The $A$-algebra $B$ induced by $\phi$ is separable if and of only $n+r\ge m$ and if $K[X]$ is generated as ideal by $P_1,\ldots,P_r$ and the maximal minors of $J_X(F,P)$.

\subsection{Remark}

1) The crucial interest of the two previous ``jacobian criterion of separability" is the bridge they build between the abstract (cf. 1.1(1)) and local (cf. 1.1(2)) notion of separability and the global and concrete ``jacobian condition".

2) The second jacobian criterion of separability is more adapted to algebric geometry than the first one. For instance, if $K$ is an algebraically closed field, $V \subset K^n$, with $n\in \N_*$, an algebraic $K$-affine set, $K[V]$ the $K$-algebra of regular functions on $V$, $F=\parenth{F_1,\ldots,F_n}:K^m\to V$ a polynomial map, with $F\in K[X]^n$, $m\in\N_*$ and $X=\parenth{X_1,\ldots,X_m}$ a system of $m$ indeterminates over $K$, $F^*:K[V]\to K[X]$ the co-morphism of $F$, this criterion says that $F$ or $F^*$ is unramified if and only if $n\ge m$ and jacobian matrix of $F$ with respect to $X$ has rank $m$  at each point of $K^m$. For instance, from this global point of view, this criterion tells us that the classical map $t\mapsto \parenth{t^3-t^2, t^2-t}$ from $K$ to a cubic $V$ with one node is unramified, even over the double point without checking this local unramification.

3) The criterion also implies that, with the notations of 2), if $F:K^m\to V$ is unramified, for any $K$-affine-sub-space $W$ of $K^m$, the restriction $F:W\to \overline{F(W)}$ of $F$ from $W$ to the Zariski closure of $F(W)$ in $K^n$ is also unramified.

4) More generally, for any $(m,n)\in \N_*^2$, any algebraically closed field  $K$, any algebraic $K$-affine set $V \subset K^m$, $W\subset K^n$, $Z\subset V$, and any unramified morphism $F:V\to W$, it can be proved that the restriction $F:Z\to\overline{F(Z)}$ of $F$ from $Z$ to the Zariski closure of $F(Z)$ in $W$ is also unramified, combining for instance proposotions (i), (ii) and (v) of \cite{AK}, Prop. 3.5, p.114.

5) In contract with this remarquable property of stability for morphisms of algebraic sets by restriction to algebraic subsets, it seems that this property is not true in general for \'etale (\emph{i.e.} flat and unramified) algebraic sets morphisms.

6) So, it is useful to have a  sure flatness criterion for unramified rings morphism. Besides a criterion deduced from Nagata's flatness criterion and another one deduced from the jacobian criterion of smoothnes in \cite{MW} (see \cite{AD1},1.1(15) and (16) for more detail), the most interesting of this kind of criterion is the following one given by A. Grothendieck in \cite{SGA1},Exposé I, Cor. 9.11.

\subsection{Theorem (Grothendieck flatness criterion for unramified ring morphisms)}

If $A$ is an equidimensional noetherian and normal ring, \emph{i.e.} a finite product of integrally closed neotherian commutative domains of the same dimension, and $B$ an equidimensional finitely generated $A$-algebra which is not the product of two rings and has the same dimension  as $A$, then any unramified ring monomorphism from $A$ to $B$ is flat.

\subsection{First isomorphism theorem (for algebras over an U.F.D., \cite{AD1}, th. 2.3)}

For any foctorial domain $A$ with caracteristic $p$ and fractions field $K$, and any ring isomorphism $\phi$ from $A$ to a commutative domain $B$, the following conditions are equivalent :

\begin{enumerate}
	\item The algebra $B$ is isomorphic to $A$.
	\item $B$ is an algebraic, monogeneous separable $A$-algebra whitout torsion and zero divisors such that $B^*=\phi(A^*)$ and $\dim_K K\ind{\otimes}{A}B\not\in p\N$
\end{enumerate}

\section{Characterization of U.F.D. isomorphisms}

\subsection{Definition}

Let $\phi:A\to B$ be a morphism of commutative rings.

We say that its co-morphism $\phi^{*}:\spec(B)\to\spec(A)$ is weakly injective if $\card \phi^{*}(J)\le 1$ for each $J$ of height one in $A$.

We also say that $\phi$ preserve primes if $\phi$ sends each prime element of $A$ to a prime one of $B$.

\subsection{U.F.D. isomorphism theorem}

For any morphism $\phi$ from a U.F.D. $A$ of characteristic $p$ to another U.F.D. $B$ the following conditions are equivalent :

\begin{enumerate}
	\item $\phi$ is an isomorphism.
	\item $\phi$ is a separable monomorphism of finite type with a weakly injective co-mophism such that $B^*=\phi(A^*)$ and the degree of the field extension $\phi(A_*)^{-1}\phi(A)\subset B_*^{-1}B$ is a finite and not divisible by $p$.
	\item $\phi$ is a flat morphism preserving primes such that $B^*=\phi(A^*)$ and $B_*^{-1}B$ is a finite and separable field extension of $\phi(A_*)^{-1}\phi(A)$, whose degree is not divisible by $p$.
\end{enumerate}

\begin{proof}
$1.\Rightarrow 2$. being trivial, let us assume 2. and consider a prime ideal $I\subset A$ of height one with a residue field $K$. According to the weak injectivity of $\phi^{*}:\spec(B)\to\spec(A)$, $J=\phi(I)B$ is primary. On the other hand, according to the separability of $B$ over $A$ and 1.1(9), (6), $B\ind{\otimes}{A}K$ is a separable, hence  a reduced $K$-algebra. furthermore, according to the injectivity of $\phi$ and 1.5, $B$ is flat over $A$, and hence the canonical map from $B/\phi(I)B = B\ind{\otimes}{A} A/I$ to $B\ind{\otimes}{A}K$ is injective, which implies that $B/\phi(I)B$ is reduced. It follows that $\phi(I)B$ is prime, and hence that $\phi$ preserves primes. Finally, according to the separability of $B$ over $A$ and 1.1(5), (6), $B_*^{-1}B$ is a finite and separable extension of $\phi(A_*)^{-1}\phi(A)$, which proves 3.

Let us now assume 3. According to the primitive element theorem and 1.1(8), there exists $a\in\phi(A_*)$ and $b\in B$ such $b$ is a primitive element of $B_*^{-1}B$ over $\phi(A_*)^{-1}\phi(A)$ and $B'=B\croch{\phi(a)^{-1}}=\phi(A)\croch{b,\phi(a)^{-1}}$ is separable over $A'=A\croch{a^{-1}}$ via the extension of $\phi$ to $A'$. So, $B'$ is a monogeneous and separable $A'$-algebra without torsion such that $\croch{B_*^{'-1}B':A_*^{'-1}A'}\not\in p\N$. Furthermore, since $A_*^{'-1}A'\ind{\otimes}{A'}B'$ is 	an integral algebra of finite type over the field $A_*^{'-1}A'$, it follows from Noether Normalisation theorem that this tensor product is a field contained in $B_*^{'-1}B'$ and containing $B'$, which means that this tensor product is $B_*^{'-1}B'$. So $\croch{A_*^{'-1}A'\ind{\otimes}{A'}B':A_*^{'-1}A'}\not\in p\N$. On the other hand, since $\phi$ preserves primes, it follows that $B^{'*}=\phi(A^{'*})$. According to the first isomorphism theorem 1.6, we obtain that $\phi:A'\to B'$ is an isomorphism, and hence that $B_*^{'-1}B'=\phi(A_*)^{-1}\phi(A)$. According to the flatness of $B$ over $\phi(A)$, it follows from the following lemma that we have the conclusion 1.
\end{proof}

\subsection{Bass Lemma (slight generalization of \cite{BA}, Cor. 1.3)}

For any and torsion free ring extension $B$ of a U.F.D. $A$ such that $A^*=A\cap B^*$, we have $B\cap\parenth{A_*^{-1}A}=A$, which means that the $A$-module $B/A$ is torsion free.

\begin{proof}
Let $b\in B\cap \parenth{A_*^{-1}A}-\{0\}$, $p$ and $q$ relatively prime elements of $A_*$ such that $b=pq^{-1}$. Since the multiplication by $p$ induces an injective endomorphism on the $A$-module $A/qA$, according to the flatness of $B$ over $A$, the multiplication by $p$ again induces an injective endomorphism on the $B$-module $B/qB$. Since $p=qb$, the last endomorphism is the zero one. It follows that $B/qB=\{0\}$, which means that $q\in A\cap B^*=A^*$. So, $b\in A$, and hence $B\cap \parenth{A_*^{-1}}=A$
\end{proof}

\section{Application to the Jacobian conjecture in any characteristic}

\subsection{Classical jacobian conjecture in any characteristic (cf. \cite{AD1}, 3.1)}

1) For $(n,d,p)\in \N_*^2\times \N$ with $p$ prime and the convention that $0$ is prime, and $K$ a commutative domain of characteristic $p$, let us call ``Classical Jacobian Conjecture in $n$ indeterminates in characteristic $p$ for polynomials endomorphisms of degree at most $d$ over $K$" the following statement denoted by $\cjc (n,p,d,K)$ : 

If $\phi$ is a monomorphism from a $K$-algebra of polynomials $A$ in $n$ indeterminates to another one $B$ such that $\deg\phi(X)\le d$ for each indeterminate $X$ of $A$, and $J(\phi)$ its jacobian matrix with respect to the indeterminates of $B$, then the following conditions are equivalent :

\begin{enumerate}
	\item $\phi$ is an isomorphism
	\item $\det J(\phi)\in K^*$ and $\croch{B_*^{-1}B:\phi(A_*)^{-1}\phi(A)}\not\in p\N$
\end{enumerate}

2) For $(n,d,p)\in \N_*^2\times\N$ with $p$ prime, we call ``Classical Jacobian Conjecture in $n$ indeterminates in characteristic $p$ for polynomial endomorphisms of degree at most $d$" the following statement denoted by $\cjc (n,p,d)$ : $\cjc(n,p,d,K)$ is true for all commutative domains $K$ of characteristic $p$.

3) For a prime $(n,p)\in\N_*\times\N$ with $p$ prime, we call ``the Classical Jacobian Conjecture in $n$ determinates in characteristic $p$" the following statement denoted $\cjc(n,p)$ : $\cjc(n,p,d)$ is true for all $d\in\N_*$.

4) For a prime $p\in\N$, we call ``the Classical Jacobian Conjecture in characteristic p" the following statement denoted by $\cjc(p)$ : $\cjc(n,p)$ is true for all $n\in\N_*$

\subsection{Na\"ive jacobian conjecture in any characteristic}

1) For $(n,p,d)\in\N_*^2\times\N$ with $p$ prime, and $K$ a commutative domain of characteristic $p$, let us call ``the Na\"ive Jacobian Conjecture in $n$ indeterminates in characteristical $p$ for polynomial endomorphisms of degree at most $d$ over $K$" the statement $\njc(n,p,d,K)$ deduced from $\cjc(n,p,d,K)$ by deleting the condition ``$\croch{B_*^{-1}B:\phi(A_*)^{-1}\phi(A)}\not\in p\N$".

2) For $(n,p,d)\in\N_*^2\times\N$ with $p$ prime, let us call ``the Na\"ive Jacobian Conjecture in $n$ indeterminates in characteristic $p$ for polynomial endomorphisms of degree at most $d$" the statement $\njc(n,p,d)$ deduced from $\cjc(n,p,d)$ by deleting the same conditions as in 1).

\subsection{Remark}

1) The Na\"ive Jacobian Conjecture in $n$ indeterminates in characteristic $p$ is trivially false according to the well-known example : $A=B=\F_p[X]$, with $X$ one indeterminate over $\F_p$, the finite field of cardinality $p$ and $\phi(X)=X-X^p$.

2) But this Na\"ivety can be easily corrected by assuming $p>d^n$, thanks to the following proposition.

\subsection{Proposition (on a reformulation of $\cjc(n,p,d,K)$ for $p>d^n$)}

For any $(n,p,d)\in\N_*^3$ with $p$ prime $p>d^n$ and $K$ a commutative domain of characteristic $p$, $\cjc(n,p,d,K)$ is equivalent to $\njc(n,p,d,K)$.

\begin{proof}
It result from the following consequence of the ``weak from of Bezot theorem" of S. Abhyankar in \cite{AB}, 12.3.1, p.272, or the ``refined Bezout theorem" of W. Fulton in \cite{FU2}, 2.3, p.10 (see also \cite{FU1}, ex. 8.4.6, p.148 and \cite{V}, p.85) :
\end{proof}

\subsection{Proposition (affine generalization of the original Bezout theorem)}

For any finite set of non constant polynomials in $n$ indeterminates over an algebraically closed field $K$, the number of isolated points of the algebraic set of $K^n$ defined by these polynomials is at most the product of their degrees.

\subsection{Corollary (a majoration of the separable geometric degree of a dominating polynomial map)}

If $B$ is an algebra of polynomials in $n>0$ indeterminates over a commutative domain $K$, $A$ the $K$-sub-algebra of $B$ generated by system of $n$ $K$-algebraically independant elements $F_1,\ldots,F_n$ of $B$, then the separable degree $\croch{B_*^{-1}B:A_*^{-1}A}_S$ of the field extension $A_*^{-1}A\subset B_*^{-1}B$ is at most $\prod_{1\le i\le n}\deg F_i$.

\begin{proof}
It follows from the previous propostition, thanks to the cardinality of the generic fiber of a dominating morphism of irreductible varieties (see for instance \cite{SP}, th. 4.1.6.)
\end{proof}

\subsection{Proposition (influence of extention of $K$ on $\cjc(n,p,d,K)$)}

If $(n,d,p)\in\N_*^2\times\N$ with $p$ prime, $K$ and $L$ commutative domains of characteristic $p$ such that $K\subset L$, then $\cjc(n,p,d,L)$ implies $\cjc(n,p,d,K)$.

\begin{proof}
Let assume $\cjc(n,p,d,L)$ and assumptions of $\cjc(n,p,d,K)$, and let $g\deg(\phi)$ be $\croch{B_*^{-1}B:\phi(A_*)^{-1}\phi(A)}$, $\phi\ind{\otimes}{K}L$ the $L$-algebra morphism from $A\ind{\otimes}{K}L$ to $B\ind{\otimes}{K}L$ induced by $\phi$. We only have to prove that condition 2. of $\cjc(n,p,d,K)$ implies that $\phi$ is an isomorphism. So, let us assume this condition. Since $g\deg(\phi)=g\deg(\phi\ind{\otimes}{K}L)$, we have $g\deg(\phi\ind{\otimes}{K}L)\not\in p\N$. On the other hand, according to 1.3 and 1.1(9), $\phi\ind{\otimes}{K}L$ is separable. So, it follows from $\cjc(n,p,d,L)$ and 1.3 that $\phi\ind{\otimes}{K}L$ is an isomorphism. According to \cite{BCW}, I, (1.1) 3,  we can conclude that $\phi$ is an isomorphism.
\end{proof}

\subsection{First reformulation mod $p$ theorem (for the Classical Jacobian Conjecture in Characteristic $0$)}

For any $(n,d)\in\N_*^2$, there exists $N(n,d)\in\N_*$ such that $\cjc(n,0,d)$ is equivalent to one of the statements ``$\cjc(n,p,d)$ for all prime $p>N(n,d)$" or ``$\njc(n,p,d)$ for all prime $p>N(n,d)$".

\begin{proof}
According to \cite{BCW}, I, (1.1)8, $\cjc(n,0,d)$ is equivalent to $\njc(n,0,\C)$. On the other hand, since $\njc(n,0,d,\C)$ is a first order proposition and since the field $\C$ is isomorphisc to the ultraproduct of the algebraic closures of prime finite fields according to the ultrafilter of the co-finite subsets of the set of non zero natural prime numbers, it follows from L\"os theorem that there exists an integer $N(n,d)\ge d^n$ such that $\njc(n,0,d,\C)$ is equivalent to ``$\njc(n,p,d,\overline{\F_p})$ for all time $p>N(n,d)$" (see for instance \cite{EK}, th. 3.1 and cor 3.2), and hence to ``$\njc(n,p,d,K)$ for all prime $p>N(n,d)$ and all algebraically closed field $K$ of characterstic $p$", according to the ``elementary equivalence" of algebraically closed fields of the same characteristic (see fo instance \cite{JL}, ch. 1, th. 1.13). So according to the proposition 3.3, $\cjc(n,0,d)$ is equivalent to ``$\cjc(n,p,d,K)$ for all prime $p>N(n,d)$ and all algebraically closed fields $K$ of characteristic $p$". Finally, the conclusion follows from the proposition 3.4 and 3.7.
\end{proof}

\subsection{Second reformulation mod $p$ theorem (for the Classical Jacobian Conjecture in characteristic $0$)}

$\cjc(0)$ is equivalent to one of the following statement.

\begin{enumerate}
	\item for any $(n,d)\in\N_*^2$, there exists $N(n,d)\in\N_*$ such that $\cjc(n,p,d,\F_p)$ is true all prime $p>N(n,d)$.
	\item for any $(n,d)\in\N_*^2$, there exists $N(n,d)\in\N_*$ such that $\njc(n,p,d,\F_p)$ is true for all prime $p>N(n,d)$ 
\end{enumerate}

\begin{proof}
The implications $1.\Leftarrow\cjc(0)\Rightarrow 2.$ follows from the first reformulation $\bmod p$ theorem, and the implications $1.\Rightarrow\cjc(0)\Leftarrow 2.$ from proposition 3.4, L\"os theorem and \cite{ES1}, th. 1.6, since the ultraproduct of prime finite fields according to the ultrafilter of the co-finite subsets of the set of prime elements of $\N_*$ is a field of characteristic $0$ isomorphic to a subfield of $\C$ (see for instance \cite{EK}, th. 5.3)
\end{proof}

\subsection{Recall (of main partial results about Jacobian Conjecture)}

1) $\cjc(2,0,100)$ is true according to \cite{MOH}.

2) With the notation of 3.2, for any $n\in\N_*$, any prime $p\neq 2$ and any commutative domain $K$ of characteristic $p$, $\njc(n,2,p,K)$ and hence $\cjc(n,2,p,K)$, is true according to \cite{WA}, th. 61 or \cite{BCW}, I, th. 2.4, and to proposition 3.7 above.

3)For any prime $p\in\N$ and any commutative domain $K$ of characteristic $p$, the statement ``$\njc(n,p,d,K)$ for any $(n,d)\in\N_*^2$'' is equivalent to  the statement ``$\cjc(n,p,3,K)$ for any $n\in\N_*$", according to \cite{BCW}, II, (2.2).

4) $\cjc(2,0)$ is true under the non necessar additional assumption that for one indeterminate $X$ of $A$, $\deg\phi(X)$ is at most the product of two primes, according to \cite{MA} and \cite{AO} corrected by \cite{NA}

5) For any $(n,d,p)\in\N_*^2\times\N$ with $p$ prime, and any commutative domain $K$, $\cjc(n,p,d,K)$ is true under the necessar additional assumption that $B$ is a monogeneous $A$-algebra, according to proposition 3.7 and \cite{FO},th. 1 for $p=0$ and \cite{AD1}, th. 2.3 for the general case.

6) For any $(n,d,p)\in\N_*^2\times\N$ with $p$ prime, and any U.F.D. $K$ of characteristic $p$, $\cjc(n,p,d,K)$ is equivalent to one of the following statements, according to \cite{AD1}, th. 3.2 :

\begin{enumerate}
	\item If $\phi$ is a monomorphism from a $K$-algebra of polynomials $A$ in $n$ indeterminates to another one $B$ such that $\deg \phi(X)\le d$ for each indeterminate $X$ of $A$, $\det J(\phi)\in K^*$ with $J(\phi)$ its jacobian matrix with respect to indeterminates of $B$, and $\croch{B_*^{-1}B:\phi(A_*)^{-1}\phi(A)}\not\in p\N$, then the $A$-sub-algebra of $B$ generated by any one of its indeterminates is an integrably closed domain.
	\item With the same assumptions as in 1., $B$ is flat over its $A$-sub-algebras generated by any one of its indeterminates.
\end{enumerate}

\subsection{Prime reformulation theorem (for the Classical Jacobian Conjecture in any characteristic)}

For any $(n,d,p)\in\N_*^2\times\N$ with $p$ prime and any U.F.D. $K$ of characteristic $p$, $\cjc(n,p,d,K)$ is equivalent to one of the following statements :

\begin{enumerate}
	\item If $\phi$ is a monomorphism from a $K$-algebra of polynomials $A$ in $n$ indeterminates to another one $B$ such that $\deg\phi(X)\le d$ for each indeterminates $X$ of $A$, $\det J(\phi)\in K^*$ with $J(\phi)$ its jacobian matrix with respect to indeterminates of $B$ and $\croch{B_*^{-1}B:\phi(A_*)^{-1}\phi(A)}\not\in p\N$, then $\phi(a)$ is a prime element of $B$ for each prime element $a$ of $A$.
	\item With the same assumptions as in 1., $\phi(a)$ is a primary element of $B$ for each prime element $a$ of $A$.
\end{enumerate}

\begin{proof}
It follows from the U.F.D. isomorphism theorem 2.2 and \cite{WA}, th. 38.
\end{proof}

\subsection{Remark}

1) If $K$ is a field of characteristic $0$ and $\phi$ a morphism from a $K$-algebra of polynomials $A$ to another $B$ with the same number of indeterminates such that $\det J(\phi)\in K^*$, with $J(\phi)$ its jacobian matrix respect to indeterminates of $B$, A. van den Essen proved in \cite{ES2}, , lemma 3.2 that $\phi(a+\lambda)$ is a prime element of $B$ for all prime elements $a$ of $A$ such that $A/aA$ is a non singular $K$-algebra and for all $\lambda$ in a co-finite subset of $K$.

2) For $p>0$, the condition $\croch{B_*^{-1}B:\phi(A_*)^{-1}\phi(A)}\not\in p\N$ in 3.11 1. is necessary for the conclusion of 3.11.1, according to the example of 3.3 (1).

\end{document}